\documentclass[11pt]{article}
\title{On Gromov's Volume Growth Theorem}
\author{Lingwei Qin}
\date{}
\usepackage[margin=1in]{geometry}
\usepackage{amsmath,amssymb,amsthm,mathtools}
\usepackage[T1]{fontenc}
\usepackage[colorlinks=true,linkcolor=blue,citecolor=blue,urlcolor=blue]{hyperref}
\allowdisplaybreaks[2]
\numberwithin{equation}{section}
\newtheorem{theorem}{Theorem}[section]
\newtheorem{proposition}[theorem]{Proposition}
\newtheorem{lemma}[theorem]{Lemma}
\newtheorem{corollary}[theorem]{Corollary}
\theoremstyle{definition}

\theoremstyle{remark}
\newtheorem{remark}[theorem]{Remark}
\newcommand{\Ric}{\operatorname{Ric}}
\newcommand{\Scal}{\operatorname{Scal}}
\newcommand{\Vol}{\operatorname{Vol}}
\newcommand{\tr}{\operatorname{tr}}
\newcommand{\Hess}{\operatorname{Hess}}
\newcommand{\supp}{\operatorname{supp}}
\newcommand{\one}{\mathbf{1}}
\newcommand{\G}{\mathsf G}
\newcommand{\E}{\mathsf E}
\newcommand{\Q}{\mathsf Q}
\newcommand{\dd}{\,d}
\newcommand{\R}{\mathbb R}
\begin{document}
\maketitle
\begin{abstract}
Let $M^n$, $n\ge3$, be complete and noncompact, with nonnegative Ricci curvature. We derive the estimate
\[
\Vol B(p,R)\le \frac{C_n}{c}(1+R)^\alpha R^{n-2}
\]
from the radial lower bound $\Scal(x)\ge c(1+d(p,x))^{-\alpha}$, where $c>0$ and $0<\alpha<2$.
The proof extends Ge's heat-covariance argument by retaining the heat occupation of the set where a chosen scalar threshold fails. Varying that threshold with the source time makes all entropy errors integrable and yields the exact exponent $n-2+\alpha$, without a logarithmic loss.
We give explicit examples showing that this exponent is optimal and also allow failure of the curvature lower bound on a set of finite volume.
\end{abstract}

\section{Introduction}\label{intro:section}
Let $(M^n,g)$, $n\ge3$, be a connected, complete, noncompact Riemannian manifold. Under nonnegative sectional curvature and the
scalar curvature bound $\Scal\ge1$, Gromov stated the uniform estimate
\[
    \sup_{p\in M}\Vol(B(p,R))\le C_nR^{n-2},\qquad R>1,
\]
and asked whether it remains valid under nonnegative Ricci curvature
\cite[p.~114, (5') and Section 2.A(b)]{Gromov}.
This estimate also implies the asymptotic bound
\[
    \limsup_{R\to\infty}\sup_{p\in M}
    \frac{\Vol(B(p,R))}{R^{n-2}}<\infty.
\]

In dimension three, Munteanu--Wang \cite{MW} proved the linear
volume-growth estimate and also studied scalar lower bounds decaying
at infinity. Chodosh--Li--Stryker \cite{CLS} subsequently gave another
proof and extended the range of decay exponents.

The conjecture in arbitrary dimension was proved independently by
Ge \cite{Ge}, Antonelli \cite{Antonelli}, and Kong--Zhu \cite{KZ}.
Ge's proof uses the heat-kernel Fisher metric and Nash entropy.
%Antonelli obtains the estimate as a special case of volume-growth bounds under positive intermediate curvature, using a quantitative Hodge obstruction and an improvement of harmonic splitting rank. Kong--Zhu use the heat-kernel pullback of the Wasserstein metric and quantitative splitting estimates, and also prove the codimension-one volume-growth bound under a uniform deficit in unit-ball volumes. 
Koirala \cite{Koirala} gives a further heat-kernel proof in the more
general setting of positive intermediate curvature and a Ricci lower
bound, recovering the scalar-curvature estimate as a special case.
%Related developments include Han's \cite{Han} heat-entropy formulation of synthetic scalar lower bounds on metric measure Dirichlet spaces.

In this paper, we use Ge's heat-covariance and entropy method to obtain two extensions. Throughout, we assume $\Ric_g\ge0$ on all of $M$.
First, we allow the uniform positive scalar lower bound to hold only outside a set of finite volume.

\begin{theorem}[Scalar curvature positive outside a finite volume set]\label{intro:compact}
Suppose that for some $\sigma>0$, $K=\{x\in M|\Scal(x) < \sigma\}$ is a measurable set with finite volume. Set $V = \Vol(K)$.
    There is $C=C(V,n,\sigma)<\infty$ such that
\begin{equation}
    \sup_{p\in M}\Vol(B(p,R))\le C R^{n-2}\qquad(R>0).
    \label{intro:compact-all-radii}
\end{equation}
Moreover,
\begin{equation}
    \limsup_{R\to\infty}\ \sup_{p\in M}
    \frac{\Vol(B(p,R))}{R^{n-2}}\le\frac{C_n}{\sigma}.
    \label{intro:compact-asymptotic}
\end{equation}
\end{theorem}

Related volume growth estimates under curvature lower bounds outside a compact set were obtained in certain dimensions by Antonelli--Xu \cite{AX} and Zhou--Zhu \cite{ZZ25}. Under nonnegative Ricci curvature, the former established linear volume growth from a positive spectral bi-Ricci lower bound in dimensions $3\le n\le5$, while the latter obtained optimal linear volume growth under a positive bi-Ricci lower bound in dimensions $3\le n\le7$.

Second, we relate a radial scalar lower bound to the volume-growth exponent.
\begin{theorem}[The radial decay estimate]\label{thm:main}
    Suppose that, for some $p\in M$, $c>0$, and $0\le\alpha\le2$,
\[
    \Scal(x)\ge c(1+d(p,x))^{-\alpha}\qquad(x\in M).
\]
Then
\begin{equation}
    \Vol(B(p,R))\le\frac{4C_n}{c}R^{n-2+\alpha}\qquad(R\ge1).
\label{eq:main}
\end{equation}
The dimensional constant can be chosen independently of $c$ and $\alpha$.

\end{theorem}

For $n=3$, the growth exponent in Theorem \ref{thm:main}, including
an outside-compact hypothesis, was already established by
Chodosh--Li--Stryker \cite[Theorem 1.1]{CLS}, building on work of
Munteanu--Wang \cite{MW}. Their theorem covers the full range
$0\le\alpha<2$, with a constant depending on the manifold and basepoint.

We sketch the main idea of the proof before giving the details.
Ge \cite{Ge} interprets the exponent in Gromov's conjecture as a loss
of two macroscopic dimensions. The heat kernel map sends a point
$x$ to the probability measure $H(x,\cdot,t)\,dV_g$.
Pulling back the Fisher information metric $g_F$ gives the normalized
tensor $\G(t)=2t\Phi_t^*g_F$.
Its trace describes the directions detected at scale $\sqrt t$, and
$2\E/t=n-\tr_g\G$ records the deficit from the Euclidean value.
The proof relates this deficit to the normalized Nash entropy $S$
through the source identities
\[
    \square_x\E=\frac12\Q,
    \qquad
    \square_x(-S)=\frac{\E}{t^2}.
\]
A curvature estimate supplies a lower bound close to two for the
heat average of $\Q$. Successive heat-potential comparisons then
transfer this bound to $\E$ and to $-S$, yielding
$S_p(T)\le-\log T+C_n$.
The entropy-to-volume inequality converts this logarithmic decay into
the growth bound $\Vol(B(p,R))\le C_nR^{n-2}$.

We apply the same method while keeping track of the region where a chosen
scalar threshold fails. Write $P_s$ for the heat semigroup and set
$B_\kappa=\{x\in M:\Scal(x)<\kappa\}$.
Proposition \ref{heat:source} gives
\[
    P_s\Q_\bullet(t)(p)
    \ge2-C_n(\kappa t)^{-1/3}-\frac{C_n}{\kappa s}
    -2P_s\one_{B_\kappa}(p),
    \qquad \kappa t\ge1.
\]
Thus the loss in the source estimate is controlled by the heat-kernel
mass of $B_\kappa$. For a fixed threshold $\kappa>0$, define
\[
    J_{\kappa,p}(T)=\int_{\kappa^{-1}}^T
    \left(\frac1t-\frac1T\right)
    P_{T-t}\one_{B_\kappa}(p)\,dt,
    \qquad T>\kappa^{-1}.
\]
The two propagation steps yield
\[
    \Vol(B(p,R))\le\frac{C_n}{\kappa}R^{n-2}
    \exp\bigl(J_{\kappa,p}(R^2)\bigr),
    \qquad R>\kappa^{-1/2}.
\]
This formula expresses the effect of the low-curvature region as an
integrated heat-kernel error. If that region lies in a set of finite volume, Gaussian estimates and
Bishop--Gromov comparison bound $J_{\kappa,p}(T)$ in terms of
$\Vol(B_\kappa)$ and $\Vol(B(p,\sqrt T))$, yielding the desired
volume estimate by a two-case argument.

Under a radial scalar lower bound, the low-curvature region can instead
be controlled by its distance from $p$. For a fixed observation radius
$R$, we take
\[
    T=R^2,\qquad a=c^{-1}(1+R)^\alpha,
    \qquad \kappa(t)=(at)^{-1/2}\quad(a\le t<T),
\]
when $a<T$. As $t$ increases, this choice moves
$B_{\kappa(t)}$ farther from $p$, while $\kappa(t)t$ increases.
The Gaussian tail bound and the source estimate then give bounded
integrated errors, with constants depending only on $n$.
Consequently,
$\Vol(B(p,R))\le C_naT^{n/2-1}$, which gives the asserted exponent.

The sharpness examples are presented in Section \ref{sharp:section}.
The product $S_a^2\times\R^{n-2}$ realizes the exponent $n-2$.
For $0<\alpha<2$, we use an explicit smooth rotationally symmetric
metric with warping function
$f(r)=r(1+r^2)^{(\alpha-2)/4}$, and take its product with
$\R^{n-3}$. This construction follows the three-dimensional
warped-product examples in \cite[Section 3]{CLS} and gives nonnegative
sectional curvature together with a positive finite limit of
$R^{-(n-2+\alpha)}\Vol(B(p,R))$.

The paper is organized as follows. Section \ref{sharp:section} constructs
the examples realizing both volume exponents.
Section \ref{heat:foundations} records the heat-kernel quantities and
analytic estimates used throughout the proof.
Section \ref{heat:source-section} establishes the source estimate with
the low-curvature occupation error.
Section \ref{heat:propagation} propagates this estimate to entropy and
volume, allowing the threshold to vary with the source time.
Sections \ref{compact:section} and \ref{decay:section} apply the common
estimate to curvature positive outside a set of finite volume and to radial
curvature decay, respectively, and include the finite-volume extensions.

\textbf{Acknowledgements} The author thanks ChatGPT for assistance with language editing and improving the clarity of the manuscript. The author also thanks Professor Jian Ge for helpful suggestions concerning theorem \ref{intro:compact} and related conclusion, in particular, for his help in shaping the final form of theorem \ref{intro:compact} .

\section{Sharpness of the volume exponents}\label{sharp:section}

The sharpness of the exponent in Theorem \ref{thm:main} can be illustrated by several model spaces. In the case $\alpha=0$, Ge \cite{Ge} already provided a sharp model, namely $S^2\times\mathbb{R}^{n-2}$, whose volume growth is of order $R^{n-2}$. In this section, we construct model spaces for the remaining cases $\alpha\neq 0$, covering both the range $0<\alpha<2$ and the endpoint $\alpha=2$. These examples demonstrate that the volume growth exponent in our estimate is optimal in the corresponding ranges of $\alpha$.Proposition \ref{sharp:decay} shows that the exponent $n-2+\alpha$ for $0 < \alpha < 2$ in Theorem \ref{thm:main} is sharp. 

\begin{proposition}\label{sharp:decay}
For every $n\ge3$ and $0<\alpha<2$, there are a smooth complete noncompact manifold with nonnegative sectional curvature, a point $p$, and a constant $c_\alpha>0$, such that
\[
    \Scal(x)\ge c_\alpha(1+d(p,x))^{-\alpha},
    \qquad
    0<\lim_{R\to\infty}\frac{\Vol(B(p,R))}{R^{n-2+\alpha}}<\infty.
\]
\end{proposition}

\begin{proof}
Put $b=\alpha/2\in(0,1)$, and define on $\R^3$
\[
    g_b=dr^2+f_b(r)^2g_{S^2},\qquad
    f_b(r)=r(1+r^2)^{(b-1)/2}.
\]
We check the metric at the origin before using the warped-product formulas.
Let $A(s)=(1+s)^{b-1}$. In Euclidean coordinates, for $r=|x|>0$,
\[
    (g_b)_{ij}=A(r^2)\delta_{ij}+\frac{1-A(r^2)}{r^2}x_ix_j.
\]
Since $A(0)=1$, the quotient $(1-A(s))/s$ extends smoothly to $s=0$.
Thus the metric extends smoothly through the origin $o$, where it equals the Euclidean inner product.
It is positive definite everywhere because its radial eigenvalue is one and its tangential eigenvalues are $A(r^2)>0$.

The radial coordinate is the distance from $o$: every curve from $o$ to radial coordinate $r$ has length at least $r$, and a radial segment has length exactly $r$.
The sets $\{r\le R\}$ are compact, so the metric is complete.

Differentiation gives
\[
    f_b'(r)=(1+br^2)(1+r^2)^{(b-3)/2}>0,
\]
\[
    f_b''(r)=(b-1)r(3+br^2)(1+r^2)^{(b-5)/2}<0
    \quad(r>0).
\]
Since $f_b'(0)=1$, one has $0<f_b'(r)\le1$.
The radial and tangential sectional curvatures are therefore nonnegative:
\[
    K_{\mathrm{rad}}=-\frac{f_b''}{f_b},
    \qquad
    K_{\mathrm{tan}}=\frac{1-(f_b')^2}{f_b^2}.
\]
For $r>0$ both are positive. Moreover,
\[
    f_b(r)=r-\frac{1-b}{2}r^3+O(r^5)
\]
shows that both curvatures tend to $3(1-b)$ at the origin.
In dimension three,
\[
    \Scal_{g_b}
    =-4\frac{f_b''}{f_b}
    +2\frac{1-(f_b')^2}{f_b^2}.
\]
It is positive everywhere, with value $18(1-b)$ at $o$.
As $r\to\infty$,
\[
f_b(r)\sim r^b,\qquad f_b'(r)\sim br^{b-1},
\qquad \frac{f_b''(r)}{f_b(r)}\sim b(b-1)r^{-2}.
\]
Since $b<1$, these formulas imply
\[
r^{2b}\Scal_{g_b}(r)\longrightarrow2.
\]
The positive continuous function
$\Scal_{g_b}(r)(1+r)^{2b}$ has a positive limit at infinity and a positive minimum on every compact interval.
Its global infimum is therefore positive.
Choose $c_\alpha>0$ below that infimum.

Now let $m=n-3\ge0$, and set
\[
(M,g)=(\R^3,g_b)\times\R^m,\qquad p=(o,0).
\]
Products preserve nonnegative sectional curvature and the Euclidean factor does not change scalar curvature.
For $x\in\R^3$ with radial coordinate $r$ and $z\in\R^m$,
\[
    d_M(p,(x,z))=\sqrt{r^2+|z|^2}\ge r.
\]
Hence the scalar lower bound persists in the full product distance.

Let $\omega_m$ denote the volume of the Euclidean unit $m$-ball, with $\omega_0=1$.
Integration over the Euclidean factor gives the exact formula
\[
    \Vol(B(p,R))=4\pi\omega_m\int_0^R
    f_b(r)^2(R^2-r^2)^{m/2}\,dr.
\]
For $m=0$, this is interpreted as the usual radial volume integral in dimension three.
Substitute $r=Rv$ to obtain
\[
    \frac{\Vol(B(p,R))}{R^{m+1+2b}}
    =4\pi\omega_m\int_0^1
    \frac{f_b(Rv)^2}{R^{2b}}(1-v^2)^{m/2}\,dv.
\]
For $r>0$,
\[
    \frac{f_b(r)}{r^b}
    =\left(\frac{r^2}{1+r^2}\right)^{(1-b)/2}\le1.
\]
Thus the integrand is bounded by the integrable function
$v^{2b}(1-v^2)^{m/2}$, and it converges to that function for every $v\in(0,1)$.
Dominated convergence yields
\[
    \lim_{R\to\infty}\frac{\Vol(B(p,R))}{R^{n-2+\alpha}}
    =4\pi\omega_{n-3}\int_0^1
    v^\alpha(1-v^2)^{(n-3)/2}\,dv.
\]
The integral is positive and finite, proving the assertion.
\end{proof}

\begin{remark}
   As $\alpha\to 2$, the model constructed in the preceding proposition degenerates to a flat model, which no longer satisfies the required scalar curvature lower bound. Thus, at the endpoint $\alpha=2$, we instead take
\[
\widetilde f(r)=\frac12 f(r)+\frac12\arctan r.
\]
    This gives a smooth model satisfying the required curvature conditions. By the same argument as above, its volume growth is of order $R^n$, showing that the exponent in the volume estimate is sharp also for $\alpha=2$.
\end{remark}
%Propositions \ref{sharp:decay} establish optimality of the exponents in Theorems \ref{compact:main} and \ref{decay:main}, respectively. They do not determine the optimal dimensional coefficients.

\section{Heat-kernel preliminaries}\label{heat:foundations}

Throughout, $(M^n,g)$ is smooth, connected, complete, noncompact, and
without boundary, with $n\ge3$ and $\Ric\ge0$.
We use $\Delta=\operatorname{div}\nabla$,
$\square_x=\partial_t-\Delta_x$.
Tensor norms are Hilbert--Schmidt norms. Constants $C_n,d_n>0$
depend only on the dimension and may change between occurrences.

Let $H(x,y,t)$ be the minimal heat kernel and set
\[
P_tf(x)=\int_MH(x,y,t)f(y)\dd V_g(y),\qquad
h=-\log H,\qquad \dd\mu_{x,t}(y)=H(x,y,t)\dd V_g(y).
\]
Stochastic completeness gives $P_t1=1$.
Derivatives in $x$ act on the pole before integration in $y$.

\subsection{Heat-kernel quantities and identities}

Following \cite{Ge}, define
\begin{align*}
    \G_x(t)&=2t\int_Md_xh\otimes d_xh\dd\mu_{x,t},\\
    D_x(t)&=\int_M\left|\Hess_xh-\frac{\G_x(t)}{2t}\right|^2\dd\mu_{x,t},\\
    \E_x(t)&=\frac t2\bigl(n-\tr_g\G_x(t)\bigr),\\
    \Q_x(t)&=4t^2D_x(t)+|g-\G_x(t)|^2+2t\langle\Ric_x,\G_x(t)\rangle,
\end{align*}
and
\[
    S_x(t)=-\int_MH(x,y,t)\log H(x,y,t)\dd V_g(y)
    -\frac n2\log(4\pi t)-\frac n2.
\]

For heat kernel map $\Phi_t=H(x,\cdot,t) d\Vol$, the pullback of Fisher information metric is 
\begin{equation*}
    \int_M d_x \log H \otimes d_x\log H d\mu_{x,t}. 
\end{equation*}

Thus $\G_x(t)=2t\Phi^*_tg_F$, the normalized Fisher pullback, and when $(M,g)$ is Euclidean space $\mathrm{R}^n$, then $\G_x(t) = g, tr_g\G_x(t) = n$. Hence, $\E_x(t)$ captures the dimensional deficit, which is then transferred to the Nash entropy $S$ through the heat operator .
The regularity arguments in \cite[Lemma B.1, Lemma 4.1 and Corollary 4.2]{Ge}
justify positive-time differentiation under these integrals. They use
$\Ric\ge0$, Gaussian bounds, and local parabolic estimates, so apply
under the present assumptions. The quantities above are finite and smooth.
We use the following consequences of
\cite[Proposition 2.4 and Proposition 4.4]{Ge} and their definitions,
\begin{equation}\label{heat:signs}
    0\le\G\le g,\qquad 0\le\E\le\frac{nt}{2},\qquad \Q\ge0,\qquad S\le0.
\end{equation}
The source identities are
\begin{equation}\label{heat:identities}
    \square_x\E=\frac12\Q,\qquad \square_x(-S)=\frac{\E}{t^2},
\end{equation}
as in \cite[Proposition 2.8 and Proposition 4.3]{Ge}.
Their proofs use the heat equation and Bochner's formula, with the
integral regularity just recorded.

The covariance gradient estimate in \cite[Proposition 3.5]{Ge} gives
\begin{equation}\label{heat:gradient}
    |\nabla\G(t)|^2\le C_n\bigl(tD(t)+t^{-1}\bigr)
    \le\frac{C_n}{t}\bigl(\Q(t)+1\bigr).
\end{equation}
The second inequality follows from $\Q\ge4t^2D$.
We retain the additive $1$, since a uniform positive lower bound for
$\Q$ is available only where the scalar threshold holds.
The integrated Li--Yau inequality in the output variable gives
\begin{equation}\label{heat:fisher}
    \int_MH(z,y,s)|\nabla_y\log H(z,y,s)|^2\dd V_g(y)
    \le\frac n{2s}\qquad(z\in M,\ s>0).
\end{equation}
See \cite{LY} and \cite[Lemma 3.11]{Ge}.

\subsection{Gaussian tails and the volume estimate}

The Li--Yau Gaussian estimate \cite{LY} gives
\begin{equation}\label{heat:gaussian}
    H(p,y,s)\le\frac{C_n}{\Vol(B(p,\sqrt s))}
    \exp\left(-\frac{d(p,y)^2}{5s}\right).
\end{equation}
Consequently, for $L\ge\sqrt s$,
\begin{equation}\label{heat:tail}
    P_s\one_{\{d(p,\cdot)>L\}}(p)
    \le C_n\exp\left(-d_n\frac{L^2}{s}\right).
\end{equation}
Indeed, by the heat kernel representation and \eqref{heat:gaussian},
\[
\begin{aligned}
    P_s\one_{\{d(p,\cdot)>L\}}(p)
    &=\int_{\{d(p,y)>L\}}H(p,y,s)\,d\Vol(y)\\
    &\leq \frac{C_n}{\Vol(B(p, \sqrt{s}))}
    \int_{\{d(p,y)>L\}}
    \exp\left(-\frac{d(p,y)^2}{5s}\right)d\Vol(y).
\end{aligned}
\]
Set $v=L/\sqrt{s}\geq1$, and decompose the integration region into annuli
\[
    A_j=\{y:jL<d(p,y)\leq(j+1)L\},
\qquad j\geq1.
\]
On each annulus,
\[
    \exp\left(-\frac{d(p,y)^2}{5s}\right)
    \leq e^{-j^2v^2/5}.
\]
Consequently, Bishop--Gromov volume comparison gives
\[
\begin{aligned}
    P_s\one_{\{d(p,\cdot)>L\}}(p)
    &\leq C_n\sum_{j=1}^{\infty}
    \frac{\Vol(B(p,(j+1)L))}{\Vol(B(p,\sqrt{s}))}
    e^{-j^2v^2/5}\\
    &\leq C_n\sum_{j=1}^{\infty}
    (j+1)^nv^ne^{-j^2v^2/5}\\
    &\leq C_n e^{-d_nv^2},
\end{aligned}
\]
where the last inequality follows by absorbing the polynomial factors into the Gaussian decay, using $v\geq1$. This proves \eqref{heat:tail}.
Also, taking $-\log$ in \eqref{heat:gaussian} and integrating against
$\mu_{p,T}$ yields
\begin{equation}\label{heat:volentropy}
    \log\frac{\Vol(B(p, \sqrt T))}{T^{n/2}}\le S_p(T)+C_n.
\end{equation}
Once $S_p(T)\le-\log(T/a)+C'$, for $0 < a <T$, 
combining this with \eqref{heat:volentropy} yields
\[
\Vol(B(p, \sqrt T))\le C'aT^{n/2-1}.
\]
For the radial estimate, we choose $T=R^2$ and
$a=c^{-1}(1+R)^\alpha$, and verify the weighted error bound
uniformly in $R$ whenever $a<T$. This gives the volume
exponent $n-2+\alpha$. 
So our goal is to obtain $S_p(T)\le-\log(T/a)+C'$. 
% 说明后面要通过 Q 的热平均得到 S 的上界，最终代入体积公式。
%\narrativegap{Explain how a lower bound for the heat-averaged source will give a volume upper bound.}

\section{From scalar curvature to a heat-averaged source}\label{heat:source-section}

The estimates in Section \ref{heat:foundations} require only nonnegative Ricci curvature. We now adapt Ge's argument to an arbitrary scalar threshold \(\kappa>0\).
On \(\{\Scal\ge\kappa\}\), a pointwise spectral estimate identifies
a lowest eigenline wherever \(\Q\) falls below \(2-\varepsilon\).
A spectral cutoff and the weighted Weitzenb\"ock estimate then control
the heat-kernel mass of this region. The resulting lower bound for
\(P_s\Q_\bullet(t)\) retains the mass of $\{\Scal<\kappa\}$ as an
explicit error term, as stated in Proposition \ref{heat:source}.

Since $\G_x(t)^\sharp$ obtained by raising the first index of $\G_x(t)$ using the metric is the self-adjoint endomorphism on $T_xM$, we can define the orthogonal projection $\Pi(x,t)$ onto the eigenspace corresponding to its smallest eigenvalue wherever the smallest eigenvalue is simple, that is, on the region$\{x\in M:\lambda_1(x,t)<\lambda_2(x,t)\}$, where $\lambda_1(x,t)$ and $\lambda_2(x,t)$ denote the two smallest eigenvalues of $\G_x(t)^\sharp$, counted with multiplicity. Its covariant derivative is determined by
\begin{equation*}
    (\nabla_X\Pi)(Y)=\nabla_X(\Pi (Y))-\Pi(\nabla_XY), \quad |\nabla\Pi|^2
    =\sum_{i,j=1}^n
    \bigl|(\nabla_{e_i}\Pi)(e_j)\bigr|^2.
    \end{equation*}
where \(\{e_i\}_{i=1}^n\) is a local orthonormal frame and these definitions are independent of the orthonormal frame.
For symmetric covariant two-tensors $A$ and $B$, we use the metric-induced inner product
\begin{equation*}
    \langle A,B\rangle
=\operatorname{tr}(A^\sharp B^\sharp),
\end{equation*}
Locally, if $\nu$ is a unit vector field spanning the image of $\Pi$, then
\begin{equation*}
    \Pi=\nu\otimes\nu^\flat,\qquad
\langle\Ric,\Pi\rangle=\Ric(\nu,\nu).
\end{equation*}

\begin{lemma}\label{heat:eigenline}
Fix $\kappa>0$, $0<\varepsilon\le1/2$, and $t>0$.
At every point where $\Scal\ge\kappa$, one has $\Q(t)\ge1$ whenever $\kappa t\ge1$.
If, in addition, $\kappa t\ge6/\varepsilon$ and $\Q(t)<2-\varepsilon$, then the eigenvalues of $\G^\sharp(t)$, ordered increasingly, satisfy
\[
    \lambda_1<\varepsilon/6,\qquad
    \lambda_2>\varepsilon/3,\qquad
    \lambda_2-\lambda_1>\varepsilon/6.
    \]
The orthogonal projection $\Pi$ onto the lowest eigenline then satisfies
\[
    \langle\Ric,\Pi\rangle>\kappa/2.
\]
\end{lemma}

\begin{proof}
Since $0\le\G\le g$, all the eigenvalues belong to $[0,1]$.
In an eigenbasis of $\G$, nonnegative Ricci curvature gives
\[
    \langle\Ric,\G\rangle
    =\sum_i\lambda_i\Ric(e_i,e_i)\ge\lambda_1\Scal.
\]
Therefore, if $\Scal\ge\kappa$ and $\kappa t\ge1$,
\[
    \Q\ge(1-\lambda_1)^2+2t\kappa\lambda_1
    \ge1+\lambda_1^2\ge1.
\]
Suppose also that $\Q<2-\varepsilon$. The nonnegative terms in $\Q$ give
\[
\kappa\lambda_1\le\langle\Ric,\G\rangle<1/t,
\]
and hence $\lambda_1<1/(\kappa t)\le\varepsilon/6$.
If $\lambda_2\le\varepsilon/3$, then
\[
    |g-\G|^2\ge(1-\varepsilon/6)^2+(1-\varepsilon/3)^2
    =2-\varepsilon+\frac{5\varepsilon^2}{36}>2-\varepsilon,
\]
which is impossible. This proves both eigenvalue assertions and the spectral gap.
Finally,
\[
    \langle\Ric,\G\rangle
    \ge\lambda_2\sum_{i=2}^n\Ric(e_i,e_i)
    =\lambda_2\bigl(\Scal-\langle\Ric,\Pi\rangle\bigr).
\]
It follows that
\[
    \Scal-\langle\Ric,\Pi\rangle<\frac3{\varepsilon t}\le\frac{\kappa}{2}.
\]
Since $\Scal\ge\kappa$, the required Ricci estimate follows.
\end{proof}

\begin{lemma}[A smooth spectral cutoff]\label{heat:cutoff}
For $0<\varepsilon\le1/2$ and $t>0$, there are a smooth function $0\le\phi\le1$ and an open set $U\subset M$ such that
\[
\phi=1\quad\hbox{if }\lambda_1\le\varepsilon/6,\ \lambda_2\ge\varepsilon/3,
\qquad
\supp\phi\subset U:=\{\lambda_2-\lambda_1>\varepsilon/48\}.
\]
On $U$, the lowest-eigenline projection $\Pi$ is smooth, and
\begin{equation}
|\nabla\phi|^2+\phi^2|\nabla\Pi|^2
\le\frac{C_n}{\varepsilon^2t}(\Q(t)+1).
\label{heat:cutoff-energy}
\end{equation}
The supported products extend smoothly by zero outside $U$.
\end{lemma}

\begin{proof}
Use the orthogonally invariant spectral cutoff constructed in
\cite[Section 3.3]{Ge}. Its support has a uniform spectral gap, and
\cite[Proposition 3.5 and Lemma 3.7]{Ge} give
\[
    |\nabla\phi|+\phi|\nabla\Pi|
    \le C_n\varepsilon^{-1}|\nabla\G|.
\]
Combining this with \eqref{heat:gradient} proves
\eqref{heat:cutoff-energy}. The gap buffer gives the smooth extension.
\end{proof}

% We also use \cite[Lemma 3.9]{Ge}. If $\Pi$ is a smooth orthogonal
% rank-one projection on an open set $U$, $w>0$ is smooth,
% $\phi$ is smooth with $\supp\phi\subset U$, and
% $\int_Mw\phi^2\dd V_g<\infty$, then
% \begin{equation}\label{heat:weighted}
% \int_Mw\phi^2\langle\Ric,\Pi\rangle\dd V_g
% \le C_n\int_Mw\left(
% |\nabla\phi|^2+\phi^2|\nabla\Pi|^2
% +\phi^2|\nabla\log w|^2\right)\dd V_g.
% \end{equation}
% The inequality is interpreted in the extended nonnegative real numbers.

\begin{lemma}[Weighted Weitzenb\"ock estimate {\cite[Lemma 3.9]{Ge}}]
\label{heat:weighted-line}
Let \(U\subset M\) be open and let \(\Pi\) be a smooth orthogonal
rank-one projection on \(TM|_U\). Suppose that \(w>0\) and \(\phi\)
are smooth functions on \(M\), with \(\supp\phi\subset U\) and
\(\int_Mw\phi^2\dd V_g<\infty\). Then
\begin{equation}\label{heat:weighted}
    \int_Mw\phi^2\langle\Ric,\Pi\rangle\dd V_g
    \le C_n\int_Mw\left(
    |\nabla\phi|^2+\phi^2|\nabla\Pi|^2
    +\phi^2|\nabla\log w|^2\right)\dd V_g.
\end{equation}
The inequality is interpreted in the extended nonnegative real numbers.
\end{lemma}

%\narrativegap{Explain how the line-field estimate controls the mass of the low-source region.}

For a scalar threshold $\kappa>0$, define
\begin{equation*}
    B_\kappa=\{x\in M|\Scal(x)<\kappa\},
\qquad q_{\kappa,z}(s)=P_s\one_{B_\kappa}(z).
\end{equation*}

The set is measurable and $0\le q_{\kappa,z}(s)\le1$.
The following estimate contains the scalar-curvature information needed for both volume results.

\begin{proposition}[The source with an occupation error]\label{heat:source}
For every $\kappa>0$, $s>0$, and $t>0$ with $\kappa t\ge1$,
\begin{equation}
P_s\Q_\bullet(t)(z)
\ge2-C_n(\kappa t)^{-1/3}-\frac{C_n}{\kappa s}
-2q_{\kappa,z}(s).
\label{heat:source-bound}
\end{equation}
\end{proposition}

\begin{proof}
Fix $0<\varepsilon\le1/2$ and initially assume $\kappa t\ge6/\varepsilon$.
Put
\[
    F=\{x \in M|\Scal(x)\ge\kappa\},\qquad A=F\cap\{\Q_x(t)<2-\varepsilon\},
\]
and abbreviate
\[
    Y=P_s\Q_\bullet(t)(z),\qquad m=P_s\one_A(z),
    \qquad q=q_{\kappa,z}(s).
\]
If $Y=\infty$, the claimed lower bound is immediate.
Suppose $Y<\infty$, and choose $\phi,\Pi$ from Lemma \ref{heat:cutoff}.
Lemma \ref{heat:eigenline} gives $\phi=1$ and $\langle\Ric,\Pi\rangle>\kappa/2$ on $A$.
With $w=H(z,\cdot,s)$, stochastic completeness gives $\int w\phi^2\le1$.
By \eqref{heat:weighted}, \eqref{heat:cutoff-energy}, and \eqref{heat:fisher},
\[
\frac{\kappa m}{2}
\le\int w\phi^2\langle\Ric,\Pi\rangle
\le\frac{C_n}{\varepsilon^2t}(Y+1)+\frac{C_n}{s}.
\]
Here nonnegative Ricci curvature was used to discard the part of the left integral outside $A$.
Hence, for dimensional $a_n,b_n>0$,
\begin{equation}
m\le\frac{a_n}{\kappa\varepsilon^2t}(Y+1)
+\frac{b_n}{\kappa s}.
\label{heat:mass}
\end{equation}
The term containing $s^{-1}$ has no factor $\varepsilon^{-2}$.

The three disjoint sets $B_\kappa,A,F\setminus A$ cover $M$.
On them $\Q$ is respectively bounded below by $0,1,2-\varepsilon$.
Consequently
\[
Y\ge m+(2-\varepsilon)(1-q-m)
=(2-\varepsilon)(1-q)-(1-\varepsilon)m.
\]
Set $\gamma=a_n/(\kappa\varepsilon^2t)$.
Inserting \eqref{heat:mass}, and using $Y\ge0$, gives
\[
(1+\gamma)Y\ge(2-\varepsilon)(1-q)-\gamma-\frac{b_n}{\kappa s}.
\]
The right-hand side divided by $1+\gamma$ equals
\[
2-\frac{\varepsilon+(2-\varepsilon)q+3\gamma+b_n/(\kappa s)}
        {1+\gamma}.
\]
The numerator is nonnegative, so
\[
Y\ge2-2q-\varepsilon-3\gamma-\frac{b_n}{\kappa s}.
\]
If $\kappa t\ge64$, take $\varepsilon=(\kappa t)^{-1/3}$.
Then $\varepsilon\le1/2$, $\kappa t\ge6/\varepsilon$, and
$\gamma=a_n(\kappa t)^{-1/3}$.
If $1\le\kappa t<64$, increase $C_n$ so that the right side of \eqref{heat:source-bound} is nonpositive, and use $Y\ge0$.
\end{proof}

\section{Entropy propagation and a common volume estimate}\label{heat:propagation}

\begin{lemma}[Minimal heat-potential comparison]\label{heat:potential}
Let $0<a<T$, and suppose that $u,f$ are nonnegative smooth functions on $M\times[a,T]$ satisfying
$(\partial_t-\Delta)u\ge f$.
Then
\[
u(p,T)\ge\int_a^T P_{T-t}f(\cdot,t)(p)\,dt.
\]
\end{lemma}

This is the minimal heat-potential comparison in \cite[Lemma 5.1]{Ge},
applied on the interval $[a,T]$. Its Dirichlet exhaustion proof applies
to every $a>0$.

\begin{lemma}[Two successive propagations]\label{heat:double}
For every $0<a<T$,
\begin{equation}
-S_p(T)\ge\frac12\int_a^T
\left(\frac1t-\frac1T\right)P_{T-t}\Q_\bullet(t)(p)\,dt.
\label{heat:double-bound}
\end{equation}
\end{lemma}

\begin{proof}
Apply Lemma \ref{heat:potential} first to
$(\partial_t-\Delta)\E=\Q/2$, and then to
$(\partial_t-\Delta)(-S)=\E/t^2$.
All signs required for the comparisons are given by \eqref{heat:signs}.
For $a<v\le T$, the first application gives
\[
\E_\bullet(v)\ge\frac12\int_a^vP_{v-t}\Q_\bullet(t)\,dt.
\]
The second application therefore gives
\[
-S_p(T)\ge
\frac12\int_a^T\int_a^v v^{-2}
P_{T-v}P_{v-t}\Q_\bullet(t)(p)\,dt\,dv.
\]
Every integrand is nonnegative. Tonelli's theorem and the semigroup law turn this into
\[
\frac12\int_a^T
\left(\int_t^T v^{-2}\,dv\right)P_{T-t}\Q_\bullet(t)(p)\,dt,
\]
which is \eqref{heat:double-bound}.
In particular, $a$ can be chosen separately for each fixed $T$.
\end{proof}
\begin{remark}
    Proposition \ref{heat:source} holds for every scalar threshold
$\kappa>0$. Thus, at each source time $t$, we may apply it with
$\kappa=\kappa(t)$, provided $\kappa(t)t\ge1$, and then integrate
the resulting inequality in $t$. The metric $g$ is fixed throughout,
and $\kappa(t)$ only selects the sublevel set
$\{\Scal<\kappa(t)\}$ and the corresponding error bounds.
Since this substitution is made after deriving the source estimate,
no derivative of $\kappa(t)$ occurs, and measurability suffices.
\end{remark}

%\narrativegap{Explain the choice of a scalar threshold at each source time and how the two applications will choose it.}

For $0<a<T$, put $w_T(t)=1/t-1/T$.
Let $\kappa:[a,T)\to(0,\infty)$ be measurable, with $\kappa(t)t\ge1$ for almost every $t$.
Define
\begin{align*}
I_1&=\int_a^T w_T(t)(\kappa(t)t)^{-1/3}\,dt,\\
I_2&=\frac1T\int_a^T\frac{dt}{t\kappa(t)},\\
J_{\kappa,p}&=\int_a^T w_T(t)
 P_{T-t}\one_{\{\Scal<\kappa(t)\}}(p)\,dt.
\end{align*}
%The measurability of the last integrand follows by applying Tonelli to the measurable function
%$(t,x)\mapsto H(p,x,T-t)\one_{\{\Scal(x)<\kappa(t)\}}$.
We obtain the following estimate.

\begin{proposition}[The common estimate]\label{heat:master}

Fix $p\in M$ and $0<a<T$. Let
$\kappa:[a,T)\to(0,\infty)$ be measurable and satisfy
$t\kappa(t)\ge1$ for almost every $t\in[a,T)$.
Then there is a constant $C_n$, depending only on $n$, such that
\begin{equation}\label{heat:master-bound}
\Vol(B(p, \sqrt T))\le C_n aT^{n/2-1}
\exp\!\left(C_n(I_1+I_2)+J_{\kappa,p}\right).
\end{equation}

\end{proposition}

\begin{proof}
If one of the upper integrals is infinite, the right-hand side is
interpreted as $+\infty$.
Assume the error integrals are finite.
Use Proposition \ref{heat:source} separately at each $t$ in \eqref{heat:double-bound}.
No derivative of $\kappa(t)$ is taken; all differential identities still concern the fixed metric and the original heat quantities.
The constant two contributes
\[
\int_a^T w_T(t)\,dt=\log(T/a)-1+a/T.
\]
For the term involving the remaining heat time,
\[
\frac{w_T(t)}{\kappa(t)(T-t)}
=\frac1{Tt\kappa(t)}.
\]
Thus
\[
-S_p(T)\ge\log(T/a)-1+a/T-C_n(I_1+I_2)-J_{\kappa,p}.
\]
Combining with \eqref{heat:volentropy} and exponentiating proves the claim.
\end{proof}

\begin{corollary}[A fixed scalar threshold]\label{heat:fixed}
Fix $\sigma>0$, put $a_\sigma=\sigma^{-1}$, and, for $T>a_\sigma$, define
\[
J_{\sigma,p}(T)=\int_{a_\sigma}^T
w_T(t)P_{T-t}\one_{\{\Scal<\sigma\}}(p)\,dt.
\]
Then
\begin{equation}
\Vol(B(p, R))\le\frac{C_n}{\sigma}R^{n-2}
             e^{J_{\sigma,p}(R^2)}
\qquad(R>\sigma^{-1/2}).
\label{heat:fixed-volume}
\end{equation}
\end{corollary}

\begin{proof}
In Proposition \ref{heat:master}, take $a=\sigma^{-1}$ and $\kappa(t)=\sigma$.
Writing $v=\sigma t$, one obtains
\[
I_1\le\int_1^\infty v^{-4/3}\,dv=3,
\qquad
I_2=\frac{\log(\sigma T)}{\sigma T}\le e^{-1}.
\]
Both are absorbed into $C_n$. The prefactor is
$aT^{n/2-1}=\sigma^{-1}R^{n-2}$.
\end{proof}

\section{Scalar curvature positive outside a set of finite volume}\label{compact:section}

\begin{proof}[Proof of Theorem \ref{intro:compact}]
The exceptional set $F=\{\Scal<\sigma\}$ is contained in $K$.
Set $A_K=\Vol(K)<\infty$ and $a_0=\sigma^{-1}$.
Integrating \eqref{heat:gaussian} over $F$ gives
\begin{equation*}
    q_{\sigma,p}(s)
    \le \frac{C_nA_K}{\Vol(B(p,\sqrt{s}))}.
\end{equation*}

Fix $T\ge2a_0$. For $a_0\le t\le T/2$, put $s=T-t$.
Then $T/2\le s\le T$, and Bishop--Gromov comparison gives
\begin{equation}\label{heat::bound}
    \frac{C_nA_K}{\Vol(B(p,\sqrt{s}))}
    \le
    \left(\frac{T}{s}\right)^{n/2}
    \frac{C_nA_K}{\Vol(B(p,\sqrt{T}))}
    \le
    \frac{2^{n/2}C_nA_K}{\Vol(B(p,\sqrt{T}))}.
\end{equation}
Splitting the integral defining $J_{\sigma,p}(T)$ at $T/2$,
using $w_T(t)\le1/t$, and applying \eqref{heat::bound} on the
first interval and $q_{\sigma,p}\le1$ on the second, we obtain
\begin{align*}
    J_{\sigma,p}(T)
    &\le
    \frac{C_nA_K}{\Vol(B(p,\sqrt{T}))}
    \int_{a_0}^{T/2}\frac{dt}{t}
    +\int_{T/2}^T\frac{dt}{t}\\
    &\le
    \frac{C_nA_K}{\Vol(B(p,\sqrt{T}))}
    \log(T/a_0)+\log2.
\end{align*}
Here and below, dimensional factors are absorbed into $C_n$.

Corollary \ref{heat:fixed}, with $T=R^2$, therefore gives
\begin{equation*}
    \Vol(B(p,R))
    \le \frac{C_n}{\sigma}R^{n-2}
    \exp\left(
    \frac{C_nA_K\log(\sigma R^2)}{\Vol(B(p,R))}
    \right)
    \qquad(R\ge\sqrt{2/\sigma}).
\end{equation*}
If
$\Vol(B(p,R))>C_nA_K\log(\sigma R^2)$,
the exponent is at most one, so
\begin{equation*}
    \Vol(B(p,R))\le\frac{C_n}{\sigma}R^{n-2}.
\end{equation*}
In the complementary case,
$\Vol(B(p,R))\le C_nA_K\log(\sigma R^2)$.
Combining these alternatives yields
\begin{equation}\label{compact:log-error}
    \Vol(B(p,R))
    \le \frac{C_n}{\sigma}R^{n-2}
    +C_nA_K\log(\sigma R^2)
    \qquad(R\ge\sqrt{2/\sigma}).
    \end{equation}
Since $n\ge3$,
\begin{equation*}
    \frac{\log(\sigma R^2)}{R^{n-2}}\longrightarrow0
    \qquad(R\to\infty).
\end{equation*}
Taking the supremum over $p$ in \eqref{compact:log-error}
therefore proves \eqref{intro:compact-asymptotic}.

For the all-radius estimate, first suppose $0<\sigma\le1$.
When $R\ge\sqrt{2/\sigma}$, we have $R\ge1$ and
\begin{equation*}
    \log(\sigma R^2)
    =\log\sigma+2\log R
    \le2R
    \le2R^{n-2}.
\end{equation*}
Thus \eqref{compact:log-error} gives
\begin{equation*}
    \Vol(B(p,R))
    \le\frac{C_n}{\sigma}(1+A_K)R^{n-2}.
\end{equation*}
For $0<R<\sqrt{2/\sigma}$, Bishop--Gromov comparison gives
\begin{equation*}
    \Vol(B(p,R))
    \le\omega_nR^n
    \le\frac{2\omega_n}{\sigma}R^{n-2}.
\end{equation*}
After increasing $C_n$, the same bound holds for every $R>0$.

If $\sigma>1$, then
$\{\Scal<1\}\subset\{\Scal<\sigma\}\subset K$.
Applying the preceding all-radius estimate with threshold one yields
\begin{equation*}
    \Vol(B(p,R))\le C_n(1+A_K)R^{n-2}
    \qquad(R>0).
\end{equation*}
Consequently,
\begin{equation*}
    \sup_{p\in M}\Vol(B(p,R))
    \le C_n(1+A_K)\max\{1,\sigma^{-1}\}R^{n-2}
    \qquad(R>0),
\end{equation*}
which proves \eqref{intro:compact-all-radii}.
\end{proof}

\begin{corollary}[compact exceptional set at a fixed threshold]\label{compact:finite}
Let $K\subset M$ be compact and suppose that, for some $\sigma>0$,
\[
    \Scal\ge\sigma\qquad\hbox{on }M\setminus K.
\] 
Set $V=\Vol(K)$. There is $C=C(V,n,\sigma)<\infty$ such that
\begin{equation}
\sup_{p\in M}\Vol(B(p,R))\le C R^{n-2}\qquad(R>0).
\label{compact:all-radii}
\end{equation}
Moreover,
\begin{equation}
    \limsup_{R\to\infty}\ \sup_{p\in M}
    \frac{\Vol(B(p,R))}{R^{n-2}}\le\frac{C_n}{\sigma}.
\label{compact:asymptotic}
\end{equation}
\end{corollary}

\begin{proof}
It follows from Theorem\ref{intro:compact} and the finiteness of volume of compact set.
\end{proof}

\section{Radial scalar curvature lower bounds}\label{decay:section}

For a constant scalar threshold, the preceding argument only needed to control the heat-kernel mass of one exceptional set.
A radial lower bound supplies a family of thresholds.
We now choose them so that all three errors in Proposition \ref{heat:master} remain bounded.

\begin{theorem}[The radial decay estimate]\label{decay:main}
Suppose that, for some $p\in M$, $c>0$, and $0<\alpha\le2$,
\[
\Scal(x)\ge c(1+d(p,x))^{-\alpha}\qquad(x\in M).
\]
Then
\begin{equation}
    \Vol(B(p,R))\le\frac{C_n}{c}(1+R)^\alpha R^{n-2}
    \qquad(R>0).
\label{decay:bound}
\end{equation}
The dimensional constant can be chosen independently of $c$ and $\alpha$.
In particular,
\[
    \Vol(B(p,R))\le\frac{4C_n}{c}R^{n-2+\alpha}\qquad(R\ge1).
\]
\end{theorem}

\begin{proof}
Fix $R>0$ and set
\[
    T=R^2,\qquad a=\frac{(1+R)^\alpha}{c}.
\]
If $a\ge T$, then Bishop--Gromov gives
\[
    \Vol(B(p,R))\le\omega_nR^n\le\omega_naR^{n-2},
\]
which is \eqref{decay:bound}.
It remains to treat $a<T$.

Set $U=T/a>1$.
For $a\le t<T$, let
\[
    u=t/a,\qquad\kappa(t)=\frac1{a\sqrt u}.
\]
The condition required in Proposition \ref{heat:master} holds because
$\kappa(t)t=\sqrt u\ge1$.
For a point $x$ with $\Scal(x)<\kappa(t)$, the scalar lower bound gives
\[
    d(p,x)>(c/\kappa(t))^{1/\alpha}-1
    =(1+R)u^{1/(2\alpha)}-1
    \ge R u^{1/(2\alpha)}.
    \]
The last inequality uses $u\ge1$.
Let $s=T-t$.
Since $R u^{1/(2\alpha)}\ge R\ge\sqrt s$, the Gaussian tail estimate \eqref{heat:tail} gives
\begin{equation}
    P_{T-t}\one_{\{\Scal<\kappa(t)\}}(p)
    \le C_n\exp(-d_nu^{1/\alpha}).
    \label{decay:occupation}
\end{equation}

We estimate each term in Proposition \ref{heat:master}.
The covariance error satisfies
\[
    I_1\le\int_a^T\frac1t u^{-1/6}\,dt
    =\int_1^Uu^{-7/6}\,du\le6.
\]
For the Fisher error,
\[
    I_2=\frac1T\int_a^T\frac{dt}{\kappa(t)t}
    =\frac1U\int_1^Uu^{-1/2}\,du
    =\frac{2(\sqrt U-1)}U\le\frac12.
\]
Indeed, writing $v=\sqrt U\ge1$, the final function is
$2(v-1)/v^2$, whose maximum is $1/2$ at $v=2$.
Finally, \eqref{decay:occupation} and $w_T(t)\le1/t$ imply
\begin{align*}
    J_{\kappa,p}
    &\le C_n\int_1^\infty e^{-d_nu^{1/\alpha}}\frac{du}{u}\\
    &=C_n\alpha\int_1^\infty e^{-d_nv}\frac{dv}{v}
\le C_n,
\end{align*}
because $0<\alpha\le2$.
Thus all three errors are bounded by dimensional constants.
Equation \eqref{heat:master-bound} now gives
\[
    \Vol(B(p,R))\le C_n aT^{n/2-1}
    =\frac{C_n}{c}(1+R)^\alpha R^{n-2}.
\]
For $R\ge1$, use $(1+R)^\alpha\le2^\alpha R^\alpha\le4R^\alpha$.
\end{proof}

\begin{proof}[Proof of Theorem\ref{thm:main}]
    For $\alpha = 0$, $\Scal \ge c$, the conclusion follows from \cite[Theorem 1.1]{Ge} by scaling.

    For $0 < \alpha \le 2 $,  apply Theorem \ref{decay:main} and use
\((1+R)^\alpha\le4R^\alpha\) for \(R\ge1\).
\end{proof}

\begin{corollary}[Radial decay with a finite-volume exceptional set]\label{decay:finite}
Suppose $p\in M$, $c>0$, $0<\alpha<2$, and $F\subset M$ is measurable with finite volume. Set $V = \Vol(F)$.
Assume
\[
    \Scal(x)\ge c(1+d(p,x))^{-\alpha}\qquad(x\in M\setminus F).
\]
Then there is $C=C(n,V,c,\alpha)<\infty$ such that
\[
    \Vol(B(p,R))\le CR^{n-2+\alpha}\qquad(R>0),
\]
and
\[
    \limsup_{R\to\infty}\frac{\Vol(B(p,R))}{R^{n-2+\alpha}}\le C_n/c.
\]
In particular, the conclusion holds when \(F\) is compact.
\end{corollary}

\begin{proof}
Set $A_F=\Vol(F)$, $T=R^2$, and $a=c^{-1}(1+R)^\alpha$.
If $T\le2a$, Bishop--Gromov gives
\begin{equation*}
    \Vol(B(p,R))\le\omega_nR^n
    \le\frac{2\omega_n}{c}(1+R)^\alpha R^{n-2}.
\end{equation*}
Suppose $T>2a$ and retain $u=t/a$ and $\kappa(t)=1/(a\sqrt u)$
from the proof of Theorem \ref{decay:main}. Then
\begin{equation*}
\{\Scal<\kappa(t)\}
    \subset F\cup\{x:d(p,x)>R u^{1/(2\alpha)}\}.
\end{equation*}
The estimates for $I_1$, $I_2$, and the Gaussian tail are unchanged.
The additional occupation error is
\begin{equation*}
    J_F(T;a)=\int_a^T w_T(t)P_{T-t}\one_F(p)\,dt.
\end{equation*}
For $a\le t\le T/2$, the Gaussian estimate and Bishop--Gromov
comparison, as in \eqref{heat::bound}, give
\begin{equation*}
    P_{T-t}\one_F(p)\le\frac{C_n A_F}{\Vol(B(p,\sqrt T))}.
\end{equation*}
Using $P_{T-t}\one_F\le1$ on the remaining interval, we obtain
\begin{equation*}
J_F(T;a)\le\frac{C_n A_F}{\Vol(B(p,R))}\log(T/a)+\log2.
\end{equation*}
Proposition \ref{heat:master} now yields
\begin{equation*}
    \Vol(B(p,R))\le\frac{C_n}{c}(1+R)^\alpha R^{n-2}
    \exp\left(\frac{C_n A_F\log(T/a)}{\Vol(B(p,R))}\right).
\end{equation*}
Distinguishing whether the exponent is at most $1$, as in the proof
of Theorem \ref{intro:compact}, gives
\begin{equation*}
    \Vol(B(p,R))\le\frac{C_n}{c}(1+R)^\alpha R^{n-2}
    +C_n A_F\log\left(\frac{cR^2}{(1+R)^\alpha}\right)
    \qquad(T>2a).
\end{equation*}
Since $0<\alpha<2$, we have $T/a\to\infty$ as $R\to\infty$.
The logarithmic term is $o(R^{n-2+\alpha})$, and
$(1+R)^\alpha/R^\alpha\to1$, proving the asserted limit superior.
The all-radius estimate follows by combining this asymptotic bound
with the Bishop--Gromov inequality
$\Vol(B(p,R))\le\omega_nR^n$ for bounded $R$.
\end{proof}

\medskip

\noindent \textsc{Lingwei Qin}

\noindent Email: 2240502179@cnu.edu.cn


\newpage

\begin{thebibliography}{9}
\bibitem[Ant26]{Antonelli}
G. Antonelli, \emph{Universal Volume Growth Bounds from Positive Intermediate Curvature},
arXiv:2608.14507v2 (2026).
\url{https://arxiv.org/abs/2608.14507}.

\bibitem[AX24]{AX}
G. Antonelli and K. Xu,
\emph{New spectral Bishop--Gromov and Bonnet--Myers theorems
and applications to isoperimetry},
arXiv:2405.08918 (2024), to appear in J. Eur. Math. Soc.
\url{https://arxiv.org/abs/2405.08918}.

\bibitem[CLS23]{CLS}
O. Chodosh, C. Li and D. Stryker,
\emph{Volume growth of $3$-manifolds with scalar curvature lower bounds},
Proceedings of the American Mathematical Society \textbf{151} (2023), 4501--4511.
\url{https://arxiv.org/abs/2207.13806}.

\bibitem[Ge26]{Ge}
J. Ge, \emph{Heat kernel geometry and Gromov's volume growth conjecture},
arXiv:2608.13553v1 (2026).
\url{https://arxiv.org/abs/2608.13553}.

\bibitem[G86]{Gromov}
M. Gromov, \emph{Large Riemannian manifolds},
in Curvature and Topology of Riemannian Manifolds (Katata, 1985),
Lecture Notes in Mathematics \textbf{1201}, Springer, Berlin, 1986,
108--121.
\url{https://doi.org/10.1007/BFb0075649}.

\bibitem[K26]{Koirala}
R. Koirala,
\emph{Volume Growth under Positive Intermediate Curvature and a Ricci Lower Bound},
arXiv:2608.17977 (2026).
\url{https://arxiv.org/abs/2608.17977}.

\bibitem[KZ26]{KZ}
B. Kong and X. Zhu, \emph{Positive Scalar Curvature and Volume Growth},
arXiv:2608.14438v2 (2026).
\url{https://arxiv.org/abs/2608.14438}.

\bibitem[LY86]{LY}
P. Li and S.-T. Yau, \emph{On the parabolic kernel of the Schr\"odinger operator},
Acta Mathematica \textbf{156} (1986), 153--201.
\url{https://doi.org/10.1007/BF02399203}.

\bibitem[MW26]{MW}
O. Munteanu and J. Wang,
\emph{Geometry of three-dimensional manifolds with positive scalar curvature},
American Journal of Mathematics \textbf{148} (2026), no.~1, 131--160.
\url{https://doi.org/10.1353/ajm.2026.a980770}.

\bibitem[ZZ25]{ZZ25}
J. Zhou and J. Zhu,
\emph{Optimal volume bound and volume growth for
Ricci-nonnegative manifolds with positive bi-Ricci curvature},
J. Reine Angew. Math. \textbf{821} (2025), 1--21.
\url{https://doi.org/10.1515/crelle-2024-0100}.
\end{thebibliography}
\end{document}